\documentclass[11pt]{amsart}
\usepackage{amsmath}
\usepackage{amsfonts}
\usepackage{amssymb}
\newtheorem{thm}{Theorem}[section]

\newtheorem{cor}[thm]{Corollary}

\newtheorem{remark}[thm]{Remark}
\newtheorem{lemma}[thm]{Lemma}
\newtheorem{prop}[thm]{Proposition}
\newtheorem{defn}[thm]{Definition}

\newtheorem{prob}[thm]{Problem}
\newcommand{\bb}[1]{\mathbb{#1}}
\newcommand{\cl}[1]{\mathcal{#1}}

\begin{document}
\title[Bimodule Projections]
{Equivariant Maps and Bimodule Projections}

\author[Vern I. Paulsen]{Vern I. Paulsen}
\address{Vern I. Paulsen: Department of Mathematics, University 
of Houston,
4800 Calhoun Road, Houston, TX 77204-3008 U.S.A.}
\email{vern@math.uh.edu, URL: http://www.math.uh.edu/$\sim$vern/}
\thanks{Research supported in part by a grant from the National 
Science
Foundation}
\begin{abstract}We construct a counterexample to Solel's\cite{So}
  conjecture that the range of any contractive, idempotent, MASA
  bimodule map on $B(\cl H)$ is necessarily a ternary subalgebra. Our
  construction reduces this problem to an analogous problem about the
  ranges of idempotent maps that are equivariant with respect to
a group action. Such maps are important to understand Hamana's theory of G-injective operator
spaces and G-injective envelopes.
\end{abstract}
\maketitle
\section{Introduction}
Let $\bb T$ denote the unit circle with arc-length measure, and let $L^2(\bb T)$ and $L^{\infty}(\bb T)$ denote the square-integrable functions and essentialy bounded functions, respectively. If we identify $L^{\infty}(\bb T) \subseteq B(L^2(\bb T))$ as the multiplication operators, then it is a maximal abelian subalgebra(MASA). We construct a unital, completely positive, idempotent $L^{\infty}(\bb T)$-bimodule map on $B(L^2(\bb T))$ whose range is not a C*-subalgebra of $B(L^2(\bb T)),$ thus providing a counterexample to a conjecture of Solel\cite{So}.

Solel\cite{So} proved that if $\cl H$ is a Hilbert space, $\cl M \subseteq B(\cl H)$ is a MASA,  and $\Phi:B(\cl H) \to B(\cl H)$ is a weak*-continuous, (completely) contractive, idempotent $\cl M$-bimodule map, then the range of $\Phi, \cl R(\Phi)$ is a {\em ternary subalgebra} of $B(\cl H),$ i.e., $A,B,C \in \cl R(\Phi)$ implies that $AB^*C \in \cl R(\Phi).$ For another proof of this fact see \cite{KP}. He also conjectured that the same result would hold even when $\Phi$ was not weak*-continuous. Our example provides a counterexample to this conjecture, since a ternary subalgebra containing the identity must be a C*-subalgebra.

As a first step in our construction, we construct a unital, completely
positive, idempotent map, $\Phi: \ell^{\infty}(\bb Z) \to
\ell^{\infty}(\bb Z),$ that is equivariant with respect to the natural
action of $\bb Z$ on $\ell^{\infty}(\bb Z)$ and whose range is not a
C*-subalgebra. The construction of this map uses some results from Hamana's theory
of G-injective envelopes\cite{H3}, where G is a discrete group acting
on all of the spaces.

In Section~2, we study $\bb Z$-equivariant projections on
$\ell^{\infty}(\bb Z)$

\section{$\bb Z$-Equivariant Projections}

In this section we take a careful look at $\ell^{\infty}(\bb Z) = C(\beta \bb Z)$ and study the $\bb Z$-equivariant
maps on this space that are also idempotent.
The identification of these two spaces comes by identifying a function  $f \in C(\beta \bb Z)$ 
with the vector $v = (f(n))_{n \in \bb Z} \in \ell^{\infty}(\bb Z).$ The action of $\bb Z$ given by
$\alpha(m)f(k) = f(k+m)$ corresponds to multiplication of the vector $v$ by $B^m$ where $B$ denotes the
backwards shift. This action also corresponds to
the unique extension of the map $k \to k+m$ to a homeomorphism of
$\beta \bb Z$ and so we shall denote this homeomorphism by $\omega \to
m \cdot \omega.$

A linear map $\Phi: \ell^{\infty}(\bb Z) \to \ell^{\infty}(\bb Z)$ is $\bb Z$-equivariant
if and only if it commutes with the backwards shift. Given such a map $\Phi$ we let
$\phi_n : \ell^{\infty}(\bb Z) \to \bb C$ denote the linear functional corresponding to the $n$th-component,
so that $\Phi(v) = (\phi_n(v)).$  

Note that $\Phi$ commutes with $B$ if and only if $\phi_n(v) = \phi_0(B^nv).$
Thus there is a one-to-one correspondence between linear functionals,
$\phi: \ell^{\infty}(\bb Z) \to \bb C$ and $\bb Z$-equivariant linear
maps,
$\Phi : \ell^{\infty}(\bb Z) \to \ell^{\infty}(\bb Z)$. We shall
denote the corresponding linear map by $\Phi = \Phi_{\phi}.$

Also, it is worth noting that $\Phi_{\phi}$ is a completely positive map if and only if $\phi$ is a positive, linear
functional, $\Phi_{\phi}$ is unital if and only if $\phi$ is unital
and $\Phi_{\phi}$ is completely contractive if and only if $\phi$ is contractive.

In this section we give some characterizations of idempotent $\bb Z$-equivariant maps. Note that
this reduces to a question of which linear functionals, $\phi$ will give rise to idempotents.
We shall construct
a positive, unital, $\bb Z$-equivariant idempotent, such that the range of $\Phi$ not
a $C^*$-subalgebra.

It is well-known that the range of such a map is completely isometrically isomorphic to a $C^*$-algebra, but what we are interested in is whether or not it is actually a
$C^*$-subalgebra. These questions are the analogues of Solel's results that the ranges of $weak^*$-continuous
masa bimodule idempotents on $B(H)$ are necessarily, $C^*$-subalgebras and ternary subalgebras, respectively,
of $B(H)$.

\begin{prop} Let $\Phi: \ell^{\infty}(\bb Z) \to \ell^{\infty}(\bb Z)$ be a $\bb Z$-equivariant linear map
and decompose $\phi_0= \phi^{ac} + \phi^{s}$ into its  weak*-continuous and singular parts.
If we define $\Phi^{ac}(v) = (\phi^{ac}(B^nv))$ and $\Phi^{s}(v) = ( \phi^{s}(B^nv)),$ then 
both these maps are $\bb Z$-equivariant, $\Phi^{ac}$ is weak*-continuous, $\Phi^{ac}(\cl C_0) \subseteq \cl C_0$,
$\Phi^{s}(\cl C_0) = 0,$
where $\cl C_0$ denotes the subspace of vectors that tend to 0, $\Phi = \Phi^{ac} + \Phi^{s},$ and this decomposition is
the unique decomposition of $\Phi$ into a weak*-continuous part and singular part.
\end{prop} 
\begin{proof} We only prove that $\Phi^{ac}(\cl C_0) \subseteq \cl C_0.$ There is a vector
$a=(a_n) \in \ell^1(\bb Z)$ such that $\phi^{ac}(v) = a \cdot v = \sum_{n \in \bb Z} a_nv_n.$
Hence, $\phi^{ac}(B^kv) = \sum_{n \in \bb Z} a_nv_{n+k} \to 0$ as $k \to \pm \infty.$
\end{proof}

Note that if we set $\hat{a}(e^{i \theta}) = \sum_{n \in \bb Z} a_ne^{in \theta}$, then since $a \in \ell^1(\bb Z)$
this series converges to give a continuous function on the circle. Identifying, $v \in \ell^{\infty}(\bb Z)$
with the formal series, $\hat{v}(e^{i \theta}) = \sum_{n \in \bb Z} v_n e^{in \theta},$
we have that $\widehat{\Phi^{ac}(v)}(e^{i \theta}) = \hat{a}(e^{i \theta}) \hat{v}(e^{i \theta}).$

\begin{prop} Let $\Phi: \ell^{\infty}(\bb Z) \to \ell^{\infty}(\bb Z)$ be a $\bb Z$-equivariant linear
map. If $\Phi$ is weak*-continuous and idempotent, then  $\Phi$ is either the identity map or 0.
\end{prop}
\begin{proof} We have that $\Phi = \Phi^{ac}$ and so $\Phi$ is given as "multiplication" by the continuous function
$\hat{a}$. It is easily checked that $\Phi \circ \Phi$ is also weak*-continuous and is given as multiplication by
$\hat{a}^2$. Since $\Phi$ is idempotent, $\hat{a}^2 = \hat{a}$ and since this function is continuous it must be either
constantly 0 or constantly 1, from which the result follows.
\end{proof}

\begin{thm} Let $\Phi: \ell^{\infty}(\bb Z) \to \ell^{\infty}(\bb Z)$ be a $\bb Z$-equivariant linear
map and let $\cl I$ denote the identity map. If $\Phi$ is idempotent, then either $\Phi = \Phi^{s}$ 
or $\Phi= \cl I - \Psi$ where $\Psi$ is singular and idempotent.
\end{thm}
\begin{proof} Write $\Phi = \Phi^{ac} + \Phi^{s}.$ Then $\Phi = \Phi \circ \Phi = \Phi^{ac} \circ \Phi^{ac} +
\Phi^{ac} \circ \Phi^{s} + \Phi^{s} \circ \Phi^{ac} + \Phi^{s} \circ \Phi^{s}$ and the first term in this sum is
easily seen to be weak*-continuous and each of the last three terms annihilate $\cl C_0$ and hence are singular.
Thus, by uniqueness of the decomposition, we have that $\Phi^{ac} = \Phi^{ac} \circ \Phi^{ac}.$
Hence, either $\Phi^{ac}$ is 0 or the identity. If $\Phi^{ac} = 0$, then $\Phi = \Phi^{s}$. If $\Phi^{ac}$ is the
identity, then equating the singular parts of the above equation yields, $\Phi^{s} = 2\Phi^s + \Phi^{s} \circ \Phi^{s}.$
Thus, $\Phi^{s} \circ \Phi^{s} = -\Phi^{s}$ and so $\Psi = -\Phi^{s}$ is singular and idempotent.
\end{proof}

Thus, we see that to construct all the idempotent maps, it is sufficient to construct all of the singular idempotents and the
singular part of an idempotent is either idempotent or the negative of an idempotent.

\begin{cor} Let $\Phi: \ell^{\infty}(\bb Z) \to \ell^{\infty}(\bb Z)$ be a $\bb Z$-equivariant linear
map and let $\cl I$ denote the identity map. If $\Phi$ is idempotent and contractive, then either $\Phi = \Phi^{s}$ 
or $\Phi= \cl I.$
\end{cor}
\begin{proof} We must show that if $\Phi = \cl I - \Psi$ with $\Psi$ idempotent and $\Phi$ is contractive,
then $\Psi = 0.$ Assume that $\Psi$ is not 0, and choose $v= (v_n)$ with $\|v\|=1$ and $\Psi(v) = v.$
Pick a component $k$ such that $|v_k| \ge 1/2$ and let $e_k$ denote the canonical basis vector that
is 1 in the k-th entry and 0 elsewhere.  Then $\| 2v_ke_k -v\| \le 1,$ 
but $\Phi(2v_ke_k - v) = (2v_ke_k - v) -(-v) = 2v_ke_k$ which has norm greater than 1.
\end{proof}

We would now like to define a {\em spectrum} for idempotent maps. To this end, for each $\lambda \in \bb T$,
where $\bb T$ denotes the unit circle in the complex plane, let $x_{\lambda} = (\lambda^n) \in \ell^{\infty}(\bb Z).$
Note that these vectors satisfy, $B(x_{\lambda}) = \lambda x_{\lambda}$ and that the eigenspace of B
corresponding to the eigenvalue $\lambda$ is one-dimensional. Hence, if $\Phi$ is any $\bb Z$-equivariant
linear map, then $\Phi(x_{\lambda})= c_{\lambda}x_{\lambda}$ for some scalar $c_{\lambda}.$
Moreover, if $\Phi$ is idempotent then $c^2_{\lambda} = c_{\lambda}$ and hence, $c_{\lambda}$ is 0 or 1.

\begin{defn} Let $\Phi: \ell^{\infty}(\bb Z) \to \ell^{\infty}(\bb Z)$ be a $\bb Z$-equivariant linear
map, then we set $\sigma(\Phi) = \{ \lambda \in \bb T: c_{\lambda} \ne 0 \}$ and we call this set the
{\bf spectrum} of $\Phi.$
\end{defn}

\begin{remark} Recall that every character on $\bb Z$ is of the form,
  $\rho_{\lambda}(n)= \lambda^n,$ for
  some $\lambda \in \bb T.$ Thus, under the identification between
  bounded functions on $\bb Z$ and vectors in $\ell^{\infty}(\bb Z)$,
  the vector $x_{\lambda}$ is just the character, $\rho_{\lambda}.$
 Note that since $\Phi$ is $\bb Z$-equivariant, the
  range of $\Phi, \cl R(\Phi)$ is a $\bb Z$-invariant subspace and
  $\lambda \in \sigma(\Phi)$ if and only if $x_{\lambda} \in \cl
  R(\Phi),$ i.e., if and only if $\rho_{\lambda} \in \cl R(\Phi).$
With these identifications, the set $\sigma(\Phi)$ is the same as the
  ``spectrum'' of the subspace, $sp(\cl R(\Phi))$ studied in the
  theory of spectral synthesis, although the latter definition is
  usually only made for weak*-closed subspaces.  See for example,
  \cite[Definition~1.4.1]{Be}.
\end{remark}

 We will show later that $\sigma(\Phi)$ is not always a closed subset of $\bb T$. The difficulty is that if $\lambda_n \to \lambda$ in $\bb T,$
then $x_{\lambda_n} \to x_{\lambda}$ only in the wk*-topology of $\ell^{\infty}(\bb Z),$ but the map $\Phi$ is generally singular.

\begin{prop} Let $\Phi: \ell^{\infty}(\bb Z) \to \ell^{\infty}(\bb Z)$ be a $\bb Z$-equivariant linear
map, then $\sigma(\Phi) = \{ \lambda : \phi_0(x_{\lambda}) \ne 0 \}$. If $\Phi$ is also,idempotent, then 
$\phi_0(x_{\lambda})$ is always either 0 or 1 and $\sigma(\Phi) = \{ \lambda : \phi_0(x_{\lambda}) = 1 \}.$ 
\end{prop}

The following result is fairly well-known. In particular, it can be deduced from Kadison's results on isometries of C*-algebras\cite{Ka}. We present a different proof that uses Choi's theory of multiplicative domains\cite{Ch} and our off-diagonalization method.

\begin{lemma} Let $\cl A$ be a $C^*$-algebra and let $\Phi: \cl A \to \cl A$ be
completely contractive. If $U_1, U_2,U_3$ are unitaries and $\Phi(U_i) = U_i,$
then $\Phi(U_1U^*_2U_3) = U_1U^*_2U_3.$
\end{lemma}
\begin{proof} By \cite{Pa}, there exist unital completely positive
  maps, $\Phi_i: \cl A \to \cl A, i=1,2$ such that the map, $\Psi: M_2(\cl
  A) \to M_2(\cl A)$ defined by $\Psi( \begin{pmatrix} a&b\\c&d
  \end{pmatrix}) = \begin{pmatrix} \Phi_1(a) & \Phi(b) \\ \Phi(c*)* &
  \Phi_2(d) \end{pmatrix}$ is completely positive.

Now consider the elements, $X_i = \begin{pmatrix} 0& U_i \\ 0 & 0 \end{pmatrix} \in M_2(\cl A), i=1,2,3.$ Since, $\Psi(X_i) = X_i, \Psi(X_i^*X_i)= X_i^*X_i, \Psi(X_iX_i^*)=X_iX_i^*,$ the elements $X_i, i=1,2,3$ belong to Choi's multiplicative domain \cite{Pa} of $\Psi.$

Consequently, $\begin{pmatrix} 0 & \Phi(U_1U_2^*U_3) \\ 0 & 0 \end{pmatrix} = \Psi(X_1X_2^*X_3) = X_1 \Psi(X_2^*) X_3 = \begin{pmatrix} 0 & U_1U_2^*U_3 \\ 0 & 0 \end{pmatrix},$ and the result follows.

\end{proof}

\begin{thm} Let $\Phi: \ell^{\infty}(\bb Z) \to \ell^{\infty}(\bb Z)$ be a $\bb Z$-equivariant contractive, idempotent
map. If $\lambda_1, \lambda_2, \lambda_3 \in \sigma(\Phi),$ then
$\lambda_1 \bar{\lambda_2} \lambda_3 \in \sigma(\Phi).$
For any $\lambda \in \sigma(\Phi),$ the set $\bar{\lambda} \cdot
\sigma(\Phi)$ is a subgroup of $\bb T.$
If, in addition, $\Phi$ is unital, then $\sigma(\Phi)$ is a subgroup of $\bb T.$
\end{thm}  
\begin{proof} The first statement is obvious from the above
  theorem and the fact that $x_{\lambda}$ is a unitary element of
  $\ell^{\infty}(\bb Z)$. To see the second claim let
  $\lambda \in \sigma(\Phi)$, and set $G= \bar{\lambda} \cdot \sigma(\Phi).$ Then
  $1 \in G$ and whenever, $\alpha= \bar{\lambda}\lambda_1, \mu=\bar{\lambda}\lambda_2 \in
  G,$ we have that $\alpha \cdot \mu=\bar{\lambda}(\lambda_1 \cdot
  \bar{\lambda} \cdot \lambda_2) \in G$ and $(\alpha)^{-1}=
  \bar{\alpha}=\bar{\lambda}(\lambda \cdot \bar{\lambda_1} \cdot
  \lambda) \in G.$
The final claim comes from choosing $\lambda= 1.$
\end{proof}

We will show later that $\sigma(\Phi)$ does not determine $\Phi.$ In
fact, we will give an example of a $\bb Z$-equivariant, unital
completely positive idempotent that is not the identity map for which
$\sigma(\Phi) = \bb T.$

Solel's proves\cite{So} that if $\cl M \subseteq B(\cl H)$ is a MASA
and $\Phi: B(\cl H) \to B(\cl H)$ is a
weak*-continuous contractive, idempotent $\cl M$-bimodule map, then
$\cl R(\Phi)$
is a
ternary subalgebra of $B(\cl H),$ i.e., $T,R,S \in \cl R(\Phi)$
implies that $TR^*S \in \cl R(\Phi).$ We return to this topic in more
detail in section 5. The following result is an analogue of Solel's result.
We let $\bb Z_n= \{ \lambda \in \bb T: \lambda^n=1 \}$ denote the
cyclic subgroups of order n, and let $\cl C_n= span \{ x_{\lambda}:
\lambda^n=1 \}$ denote the corresponding finite dimensional C*-subalgebras of
$\ell^{\infty}(\bb Z).$

\begin{cor} Let $\Phi: \ell^{\infty}(\bb Z) \to \ell^{\infty}(\bb Z)$
  be a $\bb Z$-equivariant contractive, idempotent map. If $\cl
  R(\Phi)$ is weak*-closed, then either $\sigma(\Phi)= \bb T$ and
  $\Phi$ is the identity map or there exists, $n$ and $\lambda \in \bb
  T$ such that $\sigma(\Phi) = \lambda \cdot \bb Z_n$ and $\cl
  R(\Phi)= x_{\lambda} \cdot \cl C_n.$ In all of these cases, $\cl
  R(\Phi)$ is a ternary subalgebra of $\ell^{\infty}(\bb Z).$
\end{cor}
\begin{proof}Let $M \subseteq \ell^{\infty}(\bb Z)$ be a weak*-closed,
  $\bb Z$-invariant subspace and let $sp(M)= \{ \lambda: x_{\lambda}
  \in M \}$, then $M$ is the weak*-closed span of $\{
  x_{\lambda}: \lambda \in sp(M) \}$. To see this, recall that if we identify $\ell^1(\bb Z)$ with the
  Weiner algebra, $A(\bb T)$, then the predual, $M_{\perp}$, is a norm
  closed ideal in $A(\bb T).$ If we let $k(\cdot)$ denote the kernel
  of an ideal and $h(\cdot)$ the hull of a set, then $h(E)^{\perp}=
  wk*-span \{ x_{\lambda}; \lambda \in E \}.$. Now $sp(M)=
  k(M_{\perp})$ and $M_{\perp}=h(sp(M))$ since $A(\bb T)$ is regular
  and semisimple. Hence, $M=
  (M_{\perp})^{\perp}= h(sp(M))^{\perp}= wk*-span \{ x_{\lambda}:
  \lambda \in sp(M) \}.$

Thus, if $\sigma(\Phi) = \bb T$, then $\cl R(\Phi)= \ell^{\infty}(\bb
Z)$ and so $\Phi$ is the identity map.

Note that since $\cl R(\Phi)$ is weak*-closed, $\sigma(\Phi)$ is a
closed subset of the circle. Thus, if $\lambda \in \sigma(\Phi),$ then
$\bar{\lambda} \cdot \sigma(\Phi) =G$ is a closed subgroup of $\bb T$
and hence, $G= \bb Z_n$ for some $n$ and it follows that $\cl R(\Phi)=
x_{\lambda} \cdot \cl C_n.$ 
\end{proof}

In spite of the above result we will later give an example of a $\bb
Z$-equivariant, contractive, unital, idempotent map whose range is not a C*-subalgebra.

We begin with one set of examples that is easy to describe, although as we will see their existence is problematical. Given a G-space $P$ and a closed G-invariant subset $Y \subseteq P$, a continuous, G-equivariant
function $\gamma :P \to Y$ is called a {\em G-retraction} provided that $\gamma(P)=Y$ and $\gamma(y) =y$ for all $y \in Y.$
In this case we also call $Y$ a {\em G-retract} of $P.$ Note that when $P$ is G-projective $Y$ is a G-retraction of $P$ if and only if $Y$ is also
G-projective.

Also, recall that the {\em corona set} of a discrete group G, is the set $\cl C(G)= \beta(G) \backslash G.$

\begin{prop}\label{2.11} The following are equivalent:
\begin{itemize}
\item[(i)] there exists a proper subset Y that is a $\bb Z$-retract of $\beta(\bb Z)$,
\item[(ii)] there exists an idempotent, $\bb Z$-equivariant, *-homomorphism, $\pi: C(\beta(\bb Z)) \to C(\beta(\bb Z)),$ that is not the identity map,
\item[(iii)] there exists a point $\omega \in \cl C(\bb Z),$ such that the closure
  of the orbit of $\omega$ is a $\bb Z$-retract of $\beta(\bb Z)$ that is contained in $\cl C(\bb Z).$
\end{itemize} 
Moreover, in this case, $\sigma(\pi)= \bb T.$
\end{prop}
\begin{proof} Clearly, (iii) implies (i). Assuming (i), let $\gamma: \beta(\bb Z) \to Y$ be the retraction map and set $\pi= \gamma^*,$ i.e., $\pi(f) = f \circ \gamma,$ and note that the fact that $\gamma$ is a $\bb Z$-retraction implies that $\pi$ is a $\bb Z$-equivariant, idempotent homomorphism.

Conversely, assuming (ii), there exists a continuous function, $\gamma: \beta(\bb Z) \to \beta(\bb Z)$ such that $\pi = \gamma^*,$ and the fact that $\pi$ is a $\bb Z$-equivariant, idempotent map, implies that $\gamma$ is a $\bb Z$-retraction. Thus, (i) and (ii) are equivalent.

Finally, assuming (ii), we have that $\pi = \gamma^*$ with $\gamma$ a
$\bb Z$-retraction onto some set $Y.$ Since $\pi$ is $\bb Z$-equivariant, there exists, $\rho: C(\beta(\bb Z)) \to \bb C,$ such that $\pi= \Phi_{\rho}.$ Since $\pi$ is a homomorphism, $\rho$ must be a homomorphism and hence there exists $\omega \in \beta(\bb Z)$ such that $\rho(f) = f(\omega).$ 

Since, for any $f \in C(\beta(\bb Z)),$ we have $f(m \cdot \omega) =
\pi(f)(m) = f(\gamma(m)),$ we see that $\gamma(m) = m \cdot \omega$
and hence the range of $\gamma$ must be the closure of the orbit of
$\omega$. Thus, the closure of the orbit of $\omega$ is the $\bb
Z$-retract, Y, and $\gamma$ is a $\bb Z$-retraction onto the orbit. 
Moreover, since $\pi$ is idempotent and not the identity, it must be singular and so, $\pi( c_0(\bb Z))= (0),$ but this implies that $Y \cap \bb Z$ is empty and so the closure of the orbit of $\omega$ is contained in $\cl C(\bb Z)).$

To see the final claim, note that since $x_{\lambda}$ is a unitary element of $C(\beta(\bb Z)),$ we have that $\pi(x_{\lambda}) \ne 0,$ and hence, $\pi(x_{\lambda}) = x_{\lambda},$ for every $\lambda \in \bb T.$
\end{proof}

We will now prove that such points and homomorphisms exist and consequently
 provide an example of a homomorphism such that $\pi$ is not uniquely determined by $\sigma(\pi).$

The construction of such a point can be deduced, essentially, from the theory of
{\em idempotent ultrafilters.} We are grateful to Gideon Schechtman
for introducing us to this theory. The usual proof of the existence of
idempotent ultrafilters is done for the semigroup, $\bb N.$ Since we
need to modify this to the case of $\bb Z$, we present a slightly different
version of this theory that avoids any reference to ultrafilters.
The following construction applies to any discrete group, but we shall
stick to $\bb Z$ for simplicity.

Recall that the homeomorphism $k \to k+1$ of $\bb Z$ extends to a
unique homeomorphism of $\beta(\bb Z)$ which we shall denote by
$\varphi.$ Note that $\varphi^{(n)}$ is the unique homeomorphic
extension of $k \to k+n.$ Also, note that $\varphi( \cl C(\bb Z))
\subseteq \cl C(\bb Z).$

Given $\omega \in \beta(\bb Z)$ the map $k \to \varphi^{(k)}(\omega)$
extends uniquely to a continuous function, $p_{\omega}: \beta(\bb Z)
\to \beta(\bb Z).$ Given $q \in \beta(\bb Z)$ we set $p_{\omega}(q)=
\omega * q.$ Note that if $\omega \in \cl C(\bb Z),$ then
$\varphi^{(n)}(\omega) \in \cl C(\bb Z)$ for all $n.$  Hence when
$\omega, q \in \cl C(\bb Z),$ then $\omega * q \in \cl C(\bb Z).$

\begin{prop} We have that $(\cl C(\bb Z), *)$ is a compact left-continuous,
  associative semigroup and there exists a point $\omega \in \cl C(\bb
  Z)$ such that $\omega * \omega = \omega.$
\end{prop}
\begin{proof} Left-continuous means that if $q_{\lambda} \to q$ then
  $\omega*q_{\lambda} \to \omega * q,$ which follows from the
  continuity of $p_{\omega}.$ Since $\cl C(\bb Z)$ is compact all that
  remains is to show that the product is associative, i.e., that
  $\omega_1*(\omega_2*q)= (\omega_1 * \omega_2) * q$ for all
  $\omega_1, \omega_2, q \in \cl C(\bb Z).$

Associativity is equivalent to proving that,
$p_{\omega_1}(p_{\omega_2}(q))= p_{\omega_1*\omega_2}(q).$ Since $\bb
Z$ is dense it will suffice to prove this equality for all $q=n \in
\bb Z.$ To this end choose a net of integers, $\{ m_{\alpha} \}$ that
converges to $\omega_2.$

We have that, $p_{\omega_1}(p_{\omega_2}(n))= p_{\omega_1}(
\varphi^{(n)}(\omega_2))= \lim_{\alpha} p_{\omega_1}(\varphi^{(n)}(m_{\alpha}))
= \lim_{\alpha} p_{\omega_1}(n+m_{\alpha})= \lim_{\alpha}
\varphi^{(n+m_{\alpha})}(\omega_1).$ Evaluating the right-hand side of
the above equation, yields, $p_{\omega_1 * \omega_2}(n) =
\varphi^{(n)}(\omega_1* \omega_2)=
\varphi^{(n)}(p_{\omega_1}(\omega_2)) =\lim_{\alpha}
\varphi^{(n)}(p_{\omega_1}(m_{\alpha}))= \lim_{\alpha}
\varphi^{(n)}(\varphi^{(m_{\alpha})}(\omega_1) = \lim_{\alpha}
\varphi^{(n+m_{\alpha})}(\omega_1)$ and so associativity follows.

The existence of the point $\omega$ now follows by \cite[Theorem~3.3]{Ber}, there exists an element, $\omega \in \cl
C(\bb Z)$ such that $\omega * \omega = \omega.$
\end{proof}

A point is called $\omega$ satisfying $\omega * \omega = \omega$ is called
an {\em idempotent point.}

\begin{thm}\label{2.13} Let $\omega \in \cl C(\bb Z)$ be an idempotent, then the
  map, $p_{\omega}: \beta(\bb Z) \to \cl R(p_{\omega})$ is a $\bb
  Z$-retraction onto a proper subset. Consequently, the
  map, $\pi_{\omega}:C(\beta(\bb Z)) \to C(\beta(\bb Z))$ defined by, $\pi_{\omega}(f) = f \circ p_{\omega}$ is
  a $\bb Z$-equivariant idempotent *-homomorphism onto a proper
  subalgebra with
  $\sigma(\pi_{\omega})= \bb T.$
\end{thm}
\begin{proof} Since, $\cl R(p_{\omega}) \subseteq \cl C(\bb Z)$ it is
  a proper subset. Note that $(p_{\omega} \circ p_{\omega})(q) =
  \omega *(\omega *q) = (\omega * \omega)*q= \omega * q =
  p_{\omega}(q),$ and so the map $p_{\omega}$ is idempotent.

Finally, to see that it is $\bb Z$-equivariant it is enough to show
that $p_{\omega} \circ \varphi= \varphi \circ p_{\omega}$. To this end, it is enough to consider
a dense set. We have that $p_{\omega}(\varphi(n))= p_{\omega}(n+1)=
\varphi^{(n+1)}(\omega) = \varphi \circ \varphi^{(n)}(\omega) =
\varphi(p_{\omega}(n)),$ and the rest of the proof follows from Proposition~\ref{2.11}.
\end{proof}

We now present an example of a completely positive, $\bb
Z$-equivariant projection, $\Phi$ such that $\sigma(\Phi)$ is a dense
subgroup of $\bb T$ that is {\em not} closed. A modification of this example will lead to the
counterexample to Solel's conjecture. 

The construction of this example uses Hamana's theory, \cite{H4}, of the
G-injective envelope, $I_G(\cl A)$ of a $C^*$-algebra, $\cl A,$ which
we will outline below.

Given a $C^*$-algebra, $\cl A$, equipped with an action by a discrete
group, $G$, 
Hamana\cite{H4} constructs a $G$-equivariant injective envelope,
$I_G(\cl A)$. Recall that maps between two spaces equipped with a
G-action are called {\em G-equivariant} if they satisfy, $\phi(g \cdot
a)= g \cdot \phi(a).$ The G-injective envelope
is a ``minimal'' injective $C^*$-algebra, $\cl B,$ containing $\cl A$,
that is also
equipped with an action of $G$.  To obtain it, one first shows that
$\cl A$ can always be $G$-equivariantly embedded into a $C^*$-algebra
equipped with a G-action that is {\em G-injective}, i.e., has the
property that G-equivariant completely positive maps extend.
It is easy to see, and is pointed out in Hamana\cite{H4}, that
$\ell^{\infty}(G)$ is always $G$-injective.

Once $\cl A$ is embedded into an object that is $G$-injective, one
obtains the $G$-injective envelope, by taking a minimal,
$G$-equivariant projection that fixes $\cl A$, just as in the
construction of the ordinary injective envelope. The key new
difference, is that one must restrict to maps that fix $\cl A$ and are
G-equivariant.

In \cite{HP} it is shown that if $\cl A$ is an abelian $C^*$-algebra,
then $I_G(\cl A)$ is also an abelian $C^*$-algebra. We give an ad hoc
argument of this fact for the case that we are interested in.

Let $\bb T$ denote the unit circle in the complex plane, fix an irrational number, $\theta_0,$ with $0< \theta_0 < 1$
and let $\lambda = e^{2 \pi  i\theta_0}$ so that the set $\{ \lambda^n \}_{n \in \bb Z}$ is dense in $\bb T$.
We regard $\bb T$ as a $\bb Z$-space with the action given by $n \cdot z = \lambda^n z.$
There exists a $\bb Z$-equivariant *-monomorphism $\Pi : C(\bb T) \to \ell^{\infty}(\bb Z)$ given by $\Pi(f) = (f(\lambda^n)).$
Dually, this *-monomorphism is induced by the continuous $\bb Z$-equivariant function $\gamma: \beta(\bb Z) \to \bb T$ that is given
uniquely by $\gamma(n) = \lambda^n.$ 

Since $C(\bb T)$ has been embedded into $\ell^{\infty}(\bb Z)$ in a
$\bb Z$-equivariant manner, we may obtain $I_{\bb Z}(C(\bb T))$ as the
range of a minimal $\bb Z$-equivariant idempotent map, $\phi,$ that fixes the
image of $C(\bb T).$ A priori, we only know that this range is an
operator subsystem of $\ell^{\infty}(\bb Z),$ but we can give it a
necessarily unique product via the Choi-Effros construction, i.e.,
for $\phi(a)$ and $\phi(b)$ in the range of $\phi$, we set $\phi(a)
\circ \phi(b)= \phi(ab).$ Note that since $\ell^{\infty}(\bb Z)$ is
abelian, this product will be abelian.
 
Thus, $I_{\bb Z}(C(\bb T))$ is an abelian $C^*$-algebra and if we
identify $I_{\bb Z}(C(\bb T)) = C(Y)$ then there is a homeomorphism,
$\eta: Y \to Y$, which gives the $\bb Z$-action on $Y, n \cdot y =
\eta^{(n)}(y),$ and on $C(Y)$ by $(n \cdot f)(y) = f(n \cdot y).$ The inclusion of $C(\bb T)$ into $C(Y)$ is given by a $\bb Z$-covering map
$h: Y \to \bb T$.

Any $\bb Z$-equivariant lifting of $\gamma, \Gamma : \beta(\bb Z) \to
Y$ gives a $\bb Z$-equivariant *-homomorphism of
$C(Y)$ into $\ell^{\infty}(\bb Z)$ with $\Pi(f) = \Gamma^*(f \circ h),$ that is, the restriction of $\Gamma^*$ to $C(\bb T)$ is $\Pi.$

Since $\Pi$ is a *-monomorphism,  $\Gamma^*$ must be a *-monomorphism
by the properties of the $\bb Z$-injective envelope. Hence, even though
$C(Y)$ was only assumed to be an operator subsystem of
$\ell^{\infty}(\bb Z)$, it can always be embedded as a
$C^*$-subalgebra.

The difficult problem, as we shall see shortly, is proving that it can
be embedded in such a way that it is {\em not} a $C^*$-subalgebra!

Since, $\Gamma^*$ is a *-monomorphism, $\Gamma$ must be an
onto map. Note that if $\Gamma(0) = y_0,$ then $\Gamma(n) = \eta^{(n)}(y_0),$ which we define to be $y_n$ and the range of $\Gamma$ is just the closure of the set $\{ y_n: n \in \bb Z \}.$ Moreover, the homomorphism is given by $\Gamma(f) = (f(y_n)) \in \ell^{\infty}(\bb Z).$

Since $C(Y)$ is $\bb Z$-injective there will exist a completely positive $\bb Z$-equivariant idempotent map $\Phi$ on $\ell^{\infty}(\bb Z)$ whose range is the image
of $C(Y), \Gamma^{*}(C(Y)).$ By the above results, $\Phi$ is either the identity map or singular. 

We now argue that $\Phi$ cannot be the identity map. To see this note that $\Pi(C(\bb T)) \cap \cl C_0(\bb Z) = (0).$ Hence, if we compose $\Pi$
with the quotient map into $\ell^{\infty}(\bb Z)/\cl C_0(\bb Z)$, then this composition will still be a *-monomorphism on $C(\bb T)$ and hence will also
be a *-monomorphism on $C(Y)$. Hence, the image of $C(Y)$ cannot be onto and hence $\Phi$ cannot be the identity map and thus is singular.

We claim that for this map, $-1 \notin \sigma(\Phi)$. First note that if we let $f_n(z) =z^n,$ then $\Pi(f_n) = x_{\lambda^n}$ and hence
$\lambda^n \in \sigma(\Phi)$ for all $n \in \bb Z.$ Thus,
$\sigma(\Phi)$ contains this dense subgroup of $\bb T$.

Now assume that $-1 \in \sigma(\Phi).$ Let $p_i \in \ell^{\infty}(\bb Z), i=0,1$ be the projections onto the even and odd integers, i.e., $p_i, i=0,1$ are the characteristic functions of these sets. By our assumption, $p_0= (x_{-1} + x_{+1})/2$ and hence, $p_1$ are in the range of $\Phi$. Hence, there exists 
disjoint, clopen sets $Y_i, i=0,1,$ with $Y= Y_0 \cup Y_1$ such that $p_i= \Gamma^*(\chi_{Y_i}), i=0,1.$

Thus, $p_0(n)= \chi_{Y_0}(y_n),$ and we see that, $y_n \in Y_0,$ for n even and
$y_n \in Y_1,$ for n odd.
 
Note that since $\Gamma$ is a
lifting of $\gamma$, we have that $h(y_0) =1 \in \bb T$ and since $\Gamma$ is
equivariant, $h(y_n) = h(n \cdot y_0)= \lambda^n.$ 

Since, $\theta_0$ was irrational, there exists a sequence of odd integers, $n_k$ such that, $\lambda^{n_k}$ converges to 1. Since $Y$ is compact, some subnet of $y_{n_k} \in Y_1$ will converge to a point, $z_0 \in Y_1,$ and $h(z_0)=1.$

Now let, $z_n= \eta^{(n)}(z_0),$ so that $h(z_n)= \lambda^n$ and define another $\bb Z$-equivariant lifting of $\gamma, \Gamma_1: \beta(\bb Z) \to Y$ by $\Gamma_1(n)= z_n.$
The homomorphism,$\Gamma_1^*$ also extends $\Pi$ and so it too must be a *-monomorphism and the orbit of $z_0$ must also be dense. 

Note that, since $z_0$ is a limit of odd $y_m$'s, we have that $z_n \in Y_0$ for n odd and $z_n \in Y_1$ for n even.  Thus, $\Gamma_1^*(\chi_{Y_0})= p_1.$

Finally, let $\Psi: C(Y) \to \ell^{\infty}(\bb Z)$ be defined by $\Psi= (\Gamma^* + \Gamma_1^*)/2.$ Then $\Psi$ is completely positive, $\bb Z$-equivariant, and it's restriction to $C(\bb T)$ is $\Pi$, so again by the properties of the $\bb Z$-injective envelope, $\Psi$ must be a complete order injection onto it's range.

But, $2 \Psi(\chi_{Y_0} - \chi_{Y_1}) = p_0- p_1 + p_1 - p_0 = 0,$ a contradiction.

Thus, $-1 \notin \sigma(\Phi).$

A similar argument can be used to show that no root of unity can be in $\sigma(\Phi).$

We summarize some of these results below:

\begin{thm} There exists a $\bb Z$-equivariant,  unital completely
  positive, idempotent, $\Phi: \ell^{\infty}(\bb Z) \to
  \ell^{\infty}(\bb Z)$ with $\cl R(\Phi)$ a C*-subalgebra, such that $\Phi$ is not a homomorphism, $\Phi(c_0(\bb Z))= 0$ and  $\sigma(\Phi)$ is a dense, proper subgroup of $\bb T.$
\end{thm}
\begin{proof} Let $\Phi$ be the projection onto $\Gamma(C(Y))$, as above. Then we have shown that $\Phi$ has the last two properties. But we have also seen that  the spectrum of a homomorphism must be the entire circle. Thus, $\Phi$ cannot be a homomorphism.
\end{proof}

We now present an example to show that the analogue of Solel's theorem
is not true in this setting.

\begin{thm}\label{rthm} There exists a $\bb Z$-equivariant, unital completely
  positive, idempotent, $\Phi: \ell^{\infty}(\bb Z) \to
  \ell^{\infty}(\bb Z)$ whose range is not a C*-subalgebra.
\end{thm}
\begin{proof} We retain the notation of the above discussion. Let $C(Y)= I_{\bb Z}(C(\bb T)), h:Y \to
  \bb T$, and $\eta: Y \to Y$ be as above. It is easy to see that if
  $h^{-1}(\{ 1 \})$ was a singleton, then necessarily, $h$ is
  one-to-one. But this is impossible since $C(\bb T)$ is not $\bb Z$-injective.
 So let $y_0 \ne w_0$ be points in $Y$ with $h(y_0)=h(w_0) =1$ and
  let $y_n= \eta^{(n)}(y_0), w_n=\eta^{(n)}(w_0).$

These points yield two continuous $\bb Z$-equivariant maps, $\Gamma_1, \Gamma_2: \beta(\bb
Z) \to Y$ by setting, $\Gamma_1(n) = y_n, \Gamma_2(n)= w_n $ and
corresponding $\bb Z$-equivariant *-homomorphisms, $\Gamma_i^*: C(Y)
\to \ell^{\infty}(\bb Z), i=1,2.$ Since both of these *-homomorphisms extend,
$\Pi= \gamma^*: C(\bb T) \to \ell^{\infty}(\bb Z)$ which is a
*-monomorphism, then by the properties of the injective envelope,
$\Gamma_i, i=1,2$
will both be *-monomorphisms.

Consider the unital, $\bb Z$-equivariant, completely positive map,
$\Psi= \frac{\Gamma_1 + \Gamma_2}{2}.$  Since $\Gamma_i, i=1,2$ are
both extensions of $\Pi$, we have that the restriction of $\Psi$ to
$C(\bb T)$ is a *-monomorphism and so by essentiality of the injective
envelope,$\Psi$ will be a complete order isomorphism of $C(Y)$ into
$\ell^{\infty}(\bb Z).$

If $\cl R(\Psi)$ was a $C^*$-subalgebra of $\ell^{\infty}(\bb Z),$
then by the Banach-Stone theorem $\Psi$ would be a *-isomorphism onto
its range. We now argue that $\Psi$ cannot be a *-homomorphism. 

Since, $C(Y)$ is injective, it is generated by its projections. Now
let, $p \in C(Y)$ be any projection, then $\Gamma_i^*(p) =
\chi_{E_i},i=1,2$ and $\Psi(p)= \chi_{E_3}$ must be the characteristic
functions of three sets. But since, $\chi_{E_3}= \frac{\chi_{E_1} +
  \chi_{E_2}}{2}$, examining points where both sides are 0 or 1,  it follows that $E_3= E_1 \cap E_2$ and $E_1^c =
(E_1 \cup E_2)^c$ and hence, $E_1=E_2=E_3.$ Thus, $\Psi=
\Gamma_1^*=\Gamma_2^*,$ which contradicts the choice of $\Gamma_1$ and
$\Gamma_2.$

Thus, $\cl R(\Psi)$ is not a C*-subalgebra, but since it is completely
order isomorphic to $C(Y)$ it is a $\bb Z$-injective operator
subsystem of $\ell^{\infty}(\bb Z)$ and so we may construct a $\bb
Z$-equivariant, completely positive, idempotent projection $\Phi$ of
$\ell^{\infty}(\bb Z)$ onto it.
\end{proof}

We are grateful to W. Arveson for the following argument.

Recall that $\bb T$ is the set of characters of $\bb Z$ and that for
each $\mu \in \bb T$ the corresponding character is the function,
$f_{\mu}(n)= \mu^n,$ i.e., $f_{\mu}= x_{\mu}$ under the identification
of functions with vectors. Thus, the $C^*$-subalgebra generated by the
set $\{ x_{\mu}: \mu \in \bb T \}$ is nothing more than the
$C^*$-subalgebra generated by the characters, which is the
$C^*$-algebra, $AP(\bb Z)$ of almost periodic functions on $\bb Z.$
By \cite{Be} $AP(\bb Z) = C(b \bb Z)$, where $b \bb Z$ is the Bohr
compactification of $\bb Z.$ Recall, $b \bb Z= \widetilde{\bb T_d},$
the group of characters of the abelian group, $\bb T_d$
where $\bb T_d$ denotes $\bb T$ with the discrete topology.

Recall that a topological space, X, is called {\it 0-dimensional} if
it is Hausdorff and
the clopen sets are a basis for the topology of X. By \cite{NZKF}It is
fairly easy to show that every extremally
disconnected, compact Hausdorff space is 0-dimensional. For a proof,
see the text, \cite[Proposition~10.69]{NZKF}.


If $C(b \bb Z)$ was injective, then $b \bb Z$ would be extremally disconnected and, consequently,
0-dimensional. Hence, by \cite[theorem~24.26]{HR}, the character group of $b
\bb Z$ would be a torsion group. But \cite[Theorem~26.12]{HR}, $\bb T_d$ is
the character group of $b \bb Z$ which is not a torsion group,
contradiction.

Therefore, $AP(\bb Z)$ is not an injective $C^*$-subalgebra of $\ell^{\infty}(\bb Z).$

This leads to the following observation.

\begin{thm} The $\bb Z$-injective envelope of $AP(\bb Z)$ is strictly
  contained in $\ell^{\infty}(\bb Z).$ There is a $\bb
  Z$-equivariant, unital completely positive projection,
  $\Phi: \ell^{\infty}(\bb Z) \to \ell^{\infty}(\bb Z)$ whose range is
  a C*-subalgebra that is *-isomorphic to $I_{\bb Z}(AP(\bb Z))$ that fixes $AP(\bb Z),$ annihilates $c_0(\bb Z)$ and $\sigma(\Phi)= \bb T.$
 \end{thm}
\begin{proof} Since $AP(\bb Z) \cap c_0(\bb Z) = 0,$ the quotient map into $\ell^{\infty}(\bb Z)/c_0(\bb Z)$ is $\bb Z$-equivariant and a complete isometry on $AP(\bb Z).$ Hence, it must be a complete isometry on both injective envelopes.

Thus, $I_{\bb Z}(AP(\bb Z))$ must be properly contained in $\ell^{\infty}(\bb Z)$ and we may choose a unital completely positive, $\bb Z$-equivariant projection onto it that annihilates $c_0(\bb Z).$ Because, $x_{\lambda} \in AP(\bb Z),$ we will have that $\Phi(x_{\lambda})= x_{\lambda},$ for every $\lambda \in \bb T.$
\end{proof}

Note that if $\omega \in \cl C(\bb Z)$ is an idempotent point, then
the induced idempotent, $\bb Z$-equivariant, *-homomorphism,
$\pi_{\omega}$ given by Theorem~\ref{2.13}
has $\sigma(\pi_{\omega})= \bb T$. Thus, $\pi_{\omega}$
is a projection that fixes $AP(\bb Z).$

\begin{prob} Is $I(AP(\bb Z))= I_{\bb Z}(AP(\bb Z))$?
\end{prob}
\begin{prob} If $\omega$ is an idempotent point, then is the range of
  $\pi_{\omega}$ a copy of $I_{\bb Z}(AP(\bb Z)),$ that is, is the
  range completely, isometrically isomorphic to $I_{\bb Z}(AP(\bb Z))$
  via a $\bb Z$-equivariant map that fixes $AP(\bb Z)$? This is
  equivalent to asking if $\pi_{\omega}$ is a minimal element in the
  set of all $\bb Z$-equivariant idempotent maps that fix $AP(\bb Z).$
\end{prob}

\section{MASA Bimodule Projections}

A subspace $\cl T \subseteq B(\cl K)$ is called a {\em ternary
  subalgebra} provided that, $A,B,C \in \cl T$ implies that $AB^*C \in
  \cl T.$
It is known that if $\Phi: B(\cl H) \to B(\cl H)$ is a completely
contractive, idempotent map, then the range of $\Phi, \cl R$ is
  completely isometrically isomorphic to a ternary sublagebra of
  operators on some Hilbert space.

Let $\cl M \subseteq B(\cl H)$ be a maximal abelian subalgebra(MASA).
It is also known that if $\Phi:B(\cl H) \to B(\cl H)$ is a MASA
bimodule map, then $||\Phi|| = ||\Phi||_{cb}.$

Solel\cite{So} proves that if $\Phi: B(\cl H) \to B(\cl H)$ is a
weak*-continuous, contractive, idempotent $\cl M$-bimodule map, then
the range of $\Phi$ is a ternary sublagebra of $B(\cl H).$ Thus, under
these stronger hypotheses, the completely isometric isomorphism is not
necessary.

In particular, Solel's result implies that the range of any
weak*-continuous, unital, completely positive, MASA bimodule
idempotent, must be a $C^*$-subalgebra of $B(\cl H).$

We will prove that the analogue of Solel's result is false in the that
non-weak*-continuous case.
When the MASA is discrete, then MASA bimodule maps are
known to be automatically weak*-continuous, so the main case of
interest is when the MASA is, for example, $L^{\infty}(\bb T)$
represented as multiplication operators on $ B(L^2(\bb T))$. This
subalgebra is maximal.

In this section, we show how $\bb Z$-equivariant idempotent maps on
$\ell^{\infty}(\bb Z)$ can be used to construct $L^{\infty}$-bimodule
idempotent maps on $B(L^2(\bb T)).$
The idea of the construction can be traced back to Arveson's
construction of a concrete projection of $B(L^2(\bb T))$ onto
$L^{\infty}(\bb T).$ 

Let $z^n=e^{in \theta}, n \in \bb Z,$ denote the
standard orthonormal basis for $L^2(\bb T).$ This basis defines a
Hilbert space isomorphism between $L^2(\bb T)$ and $\ell^2(\bb Z).$
We identify bounded operators on $\ell^2(\bb Z)$ with the infinite
matrices, $(a_{i,j}), i,j \in \bb Z$ and
using this isomorphism, the multiplication operators for a function $f$ is identified
with the bounded, {\em Laurent matrix}, $(a_{i,j})$ where, $a_{i,j}=
\hat{f}(i-j),$ the Fourier coefficient. In particular, the operator of
multiplication by $z$ corresponds to the bilateral shift operator, $B$. 

We further identify, $\ell^{\infty}(\bb Z)$ with the bounded, diagonal
operators, $\cl D \subseteq B(\ell^2(\bb Z)).$ Note that under this
identification, the action of $\bb Z$ on $\ell^{\infty}(\bb Z)$ induced
translation is implemented by conjugation by B. We define $\alpha(n): \cl
D \to \cl D$ by $\alpha(n)(D) = B^nDB^{-n}.$ Thus, a map $\Phi: \cl D
\to \cl D$ is $\bb Z$-equivariant if and only if $\Phi(B^nDB^{-n}) =
B^n \Phi(D) B^{-n},$  for all D and all $n$, which is if and only if
$\Phi(BDB^{-1})= B \Phi(D) B^{-1},$ for all D.

We define, $E: B(\ell^2(\bb Z)) \to \cl D,$ by letting $E((a_{i,j}))$
  be the diagonal operator, with diagonal entries, $a_{i,i}.$ It is
  well-known that, $E$ is a completely positive projection from $B(\ell^2(\bb
  Z))$ onto $\cl D.$ Given $A \in B(\ell^2(\bb Z))$ we set
  $\hat{A}(n)= E(AB^{-n}) \in \cl D$ and we call, $\sum_{n \in \bb Z}
  \hat{A}(n)B^n$ the {\em formal Fourier series} for A.

We should remark that just as with $L^{\infty}$-functions, the formal
Fourier series uniquely determines A, but does not need to converge to
A in any reasonable sense. In fact, it need not converge to A in even
the weak operator topology.

The key fact, whose proof we defer until later, is the following:

\begin{thm} Let $\Phi: \cl D \to \cl D$ be a $\bb Z$-equivariant,
  unital completely positive map. Then there is a well-defined unital
  completely positive, $L^{\infty}$-bimodule map, $\Gamma:
  B(\ell^2(\bb Z)) \to B(\ell^2(\bb Z))$ satisfying, $\Gamma( \sum_{n
  \in \bb Z} D_n B^n) = \sum_{n \in \bb Z} \Phi(D_n) B^n.$
Moreover, if $\Phi$ is idempotent, then $\Gamma$ is idempotent.
\end{thm}

\begin{cor} There exists a unital, completely positive, idempotent
  $L^{\infty}(\bb T)$-bimodule map, $\Gamma: B(L^2(\bb T)) \to
  B(L^2(\bb T))$ whose range is not a C*-subalgebra.
\end{cor}
\begin{proof} Let $\Phi: \cl D \to \cl D$ be the $\bb Z$-equivariant,
  unital completely positive, idempotent map given by \ref{rthm}
  whose range is not a C*-subalgebra and let $\Gamma: B(\ell^2(\bb Z))
  \to B(\ell^2(\bb Z))$ be the map given by the above theorem.
If $\cl R(\Gamma)$ was a C*-subalgebra, then $\cl R(\Gamma) \cap \cl D
  = \cl R(\Phi)$ would also be a C*-subalgebra.
Hence, $\cl R(\Gamma)$ is not a C*-subalgebra. The proof is completed
  by making the identification of $\ell^2(\bb Z)$ with $L^2(\bb T)$
  which carries the Laurent matrices to the multiplication operators.
\end{proof}

Note that when $\Phi$ is singular, $\Gamma$ will also not be weak*-continuous.
Thus, using the example from the previous section, we see that there exists
a unital, completely positive, idempotent, $L^{\infty}$-bimodule map,
$\Gamma$, such that, not only is $\Gamma(K)=0$ for every compact
operator $K$, but $\Gamma(A)= 0$ whenever, $\hat{A}(n) \in c_0(\bb Z)$
for all n.

\begin{prob} Does there exist a unital, completely positive, idempotent, $\bb
  Z$-equivariant map, $\Phi: \cl D \to \cl D,$ whose range is a
  C*-subalgebra, but such that the range of $\Gamma: B(\ell^2(\bb Z)) \to
  B(\ell^2(\bb Z))$ is not a C*-subalgebra?
\end{prob}

Before proving the above theorem, we will need a few results about
cross-products. Recall that in general we can form two crossed
products, a full and reduced crossed product, but when the group is amenable,
then these crossed-products agree \cite[Theorem~7.7.7]{Pe}.
Since $\bb Z$ is amenable, we let $\bb Z \times_{\alpha} \cl D$ denote
this C*-algebra. A dense set of elements of the crossed product is
given by finitely supported functions, $f: \bb Z \to \cl D.$ Given any pair consisting of a *-homomorphism, $\pi:
\cl D \to B(\cl H)$ and a unitary $U \in B(\cl H)$ such that $U^n
\pi(D) U^{-n}= \pi(\alpha(n)(D))$(such a pair is called a {\em
  covariant representation}), there exists a *-homomorphism, $\Pi: \bb
Z \times_{\alpha} \cl D \to B(\cl H)$ satisfying, $\Pi(f) = \sum_n f(n)U^n.$

\begin{lemma} Let $\pi: \cl D \to B(\ell^2(\bb Z))$ be the identity
  inclusion and let $B \in B(\ell^2(\bb Z))$ denote the bilateral
  shift. Then these are a covariant pair and the map, $\Pi: \bb Z
  \times_{\alpha} \cl D \to B(\ell^2(\bb Z))$ is a *-monomorphism.
\end{lemma}
\begin{proof} By \cite[Theorem~7.7.5]{Pe}, if $\lambda: \bb Z \to
  B(\ell^2(\bb Z))$ denotes the
  left regular representation, then $\tilde{\pi}= \pi \otimes id:
  \cl D \to B(\ell^2(\bb Z) \otimes \ell^2(\bb Z))$ and
  $\tilde{\lambda}= id \otimes \lambda: \bb Z \to B(\ell^2(\bb Z)
  \otimes \ell^2(\bb Z))$ are a covariant pair and the induced
  representation, $\pi \times \lambda: \bb Z \times_{\alpha} \cl D \to
  B(\ell^2(\bb Z) \otimes \ell^2(\bb Z))$ is faithful, i.e., a
  *-monomorphism.

Let $\cl H_n = e_n \otimes \ell^2(\bb Z)$ so that $\ell^2(\bb Z)
\otimes \ell^2(\bb Z) = \sum_{n \in \bb Z} \oplus \cl H_n.$ Each of
these subspaces is a reducing subspace for the image of $\bb Z
\times_{\alpha} \cl D$ and the restrictions to any pair of them are
unitarily equivalent. Moreover, the restriction to $\cl H_0$
is $\Pi$. Hence, $\Pi$ must be a *-monomorphism.
\end{proof}

Thus, we may identify $\bb Z \times_{\alpha} \cl D$ with the norm
closure in $B(\ell^2(\bb Z))$ of the operators that are finite sums of
the form, $\sum D_n B^n,$ with $D_n \in \cl D.$

\begin{thm} Let G be a discrete group, let $A$ be a C*-algebra and let
  $\alpha:G \to Aut(A)$ be a group action. If $\Phi:A \to B(\cl H)$ is
  a completely positive map, $\rho:G \to B(\cl H)$ is a unitary
  representation, such that $\rho(g)\Phi(a) \rho(g^{-1})=
  \Phi(\alpha(g)(a)),$ i.e., a covariant pair, then there is a
  completely positive map, $\rho \times \Phi: G \times_{\alpha} A \to
  B(\cl H),$ satisfying, $\rho \times \Phi(f)= \sum_{g \in G}
  \Phi(f(g)) \rho(g),$ for any finitely supported function, $f: G \to
  A.$
\end{thm}
\begin{proof} This is a restatement of the covariant version of
  Stinespring's theorem \cite[Theorem~2.1]{Pa1}. Let $\pi:A \to B(\cl
  K), \tilde{\rho}:G \to B(\cl K)$ and $V: \cl H \to \cl K$ be the
  covariant pair that dilates $\Phi, \rho,$ then $\rho \times \Phi(f)
  = V^* \tilde{\rho} \times \pi(f) V$ and the result follows.
\end{proof}

We now turn our attention to the proof of Theorem~5.1. By the above results,
 given any $\Phi: \cl D \to \cl D$ that is completely positive and $\bb Z$-equivariant,
we have a well-defined completely positive map, $\Gamma,$ satisfying
 for any finite sums,
 $\Gamma(\sum D_nB^n)= \sum \Phi(D_n)B^n$ whose domain is the norm
 closure of such finite sums and we wish to extend it to all of
 $B(\ell^2(\bb Z))$.

Consider for any $0 \le r <1$, the matrix, $P_r= (r^{|i-j|}) =
  (I-rB)^{-1} + (I- rB^*)^{-1} -I \ge 0,$ since it is the ``operator
  Poisson kernel'' \cite[Exercise~2.15]{Pa}.
Thus, the corresponding Schur product map, $S_r: B(\ell^2(\bb Z)) \to
  B(\ell^2(\bb Z))$ given by $S_r((a_{i,j}) = (r^{|i-j|}a_{i,j})$ is
  completely positive and unital. Writing $A= ( a_{i,j}) \sim \sum
  D_nB^n$ in its formal Fourier series, we see that $A_r \equiv S_r(A) = \sum
  r^{|n|}D_nB^n$ where in the latter case we have absolute norm
  convergence  of the partial sums. Hence, for any $A, S_r(A) \in \bb
  Z \times_{\alpha} \cl D.$

This shows that, $\Gamma(S_r(A))= \sum r^{|n|}\Phi(D_n)B^n.$

Note that, for any $A \in B(\ell^2(\bb Z))$ we have that $A \ge 0$ if
and only if $A_r \ge 0$ for all $0 \le r <1,$ and $\|A\| = \sup_{0 \le
  r < 1} \|A_r\|.$ Finally, given any formal matrix, $A=(a_{i,j})$ it
is easily checked that $A$ defines a bounded operator if and only if
$\sup_{0 \le r < 1} \|A_r\|$ is finite and that this supremum equals
the norm.

Hence, if $A \sim \sum D_nB^n$ is bounded, then $\|\sum
r^{|n|}\Phi(D_n)B^n\| \le \|\Phi\| \|A\|$ for all $0 \le r < 1$ and hence, $\sum
\Phi(D_n)B^n$ is bounded. This shows that we may extend $\Gamma$ to
all of $B(\ell^2(\bb Z))$ and the norm of the extended map will be at
most $\|\Phi\|.$ Using the positivity properties of $A_r$ we see that
the extended map $\Gamma$ will be positive. Complete positivity of the
extended follows similarly.

Finally, since in Theorem~5.1, we are assuming that $\Phi$ is unital,
$\Gamma$ will fix every Laurent matrix, which is a C*-subalgebra. By Choi's theory of
multiplicative domains \cite{Ch}(see also \cite{Pa}), we have that
$\Gamma$ is a bimodule map over the Laurent matrices.

This completes the proof of Theorem~5.1.

\end{document}